\documentclass[reqno]{amsart}
\usepackage{amssymb, amsmath, color, amsfonts, amsthm,verbatim, array}
\usepackage{graphicx, graphics, epsf, epsfig, quiver, enumitem, lineno}

\theoremstyle{plain}
\newtheorem{theorem}{Theorem}
\newtheorem{corollary}{Corollary}

\newtheorem{proposition}{Proposition}

\theoremstyle{definition}

\theoremstyle{remark}

\DeclareMathOperator{\im}{im}
\DeclareMathOperator{\coker}{coker}
\DeclareMathOperator{\ord}{ord}
\DeclareMathOperator{\gl}{GL}

\begin{document}

\title{Determining Compatible Groups\\ in $\Gamma$-Zonal Labelings}
\author{Andrew Bowling}
\address{Mathematics Department\\
         Denison University\\
         Granville, OH 43023\\
         USA}
\email{bowlinga@denison.edu}
\author{Richard M. Low} 
\address{Department of Mathematics and Statistics\\
         San Jose State University\\
         San Jose, CA 95192\\
         USA}
\email{richard.low@sjsu.edu}
\begin{abstract}



\date{November 16, 2026. Version 0.3}

Let $G$ be a plane graph with vertex, edge, and region sets $V(G),$ $E(G),$ and $F(G)$ respectively, and let $\Gamma$ be a nontrivial abelian group with identity element $0$. A $\Gamma$-zonal labeling of a plane graph $G$ is a labeling $\ell: V(G)\rightarrow \Gamma\setminus\{0\}$ such that for every region $f\in F(G)$, the sum of labels on the boundary of $f$ is $0$ in $\Gamma$. We say a $\Gamma$-zonal labeling is generative if the labels used generate $\Gamma$, and strong if the additive order of each label matches the degree of the corresponding vertex. Given a graph $G$ or a family of graphs $\mathcal{G}$, we explore restrictions on the groups $\Gamma$ for which there can exist a $\Gamma$-zonal labeling. Here we specifically determine restrictions on compatible groups when $\mathcal{G}$ is the family of $2$-connected multigraphs or $2$-connected simple graphs. We also demonstrate how the Smith normal form decomposition of the face incidence matrix of a plane graph $G$ can help characterize the groups $\Gamma$ for which there exists a generative $\Gamma$-zonal labeling. 
\end{abstract}
\maketitle

\section{Introduction}

 A graph is \textit{planar} if it can be drawn in the plane (or on the sphere) in such a way that no two edges cross. Any such embedding is called a \textit{plane graph.} Let $G$ be a plane graph with vertex, edge, and region sets $V(G), E(G)$ and $F(G)$, respectively. A \textit{zonal labeling} of a plane graph $G$ is a labeling $\ell: V(G)\rightarrow \mathbb{Z}_3\setminus \{0\}$ such that for every region $f\in F(G)$ with boundary $B_f$, $\sum_{v\in V(B_f)}\ell(v)=0\in \mathbb{Z}_3$. 
 
 Zonal labelings were developed by Australian physicist Cooroo Egan in 2014 \cite{CharEganZhang}. The primary motivation for developing these labelings was to provide a new perspective on studying the Four Color Theorem. However, there are many interesting results on zonal labelings that are not directly connected to this famous theorem \cite{Bowling, BowZhang, CharEganZhangRevis}. Recently, the concepts of zonality have been extended in several directions. In \cite{Bowling1}, the definition of zonal labelings is extended to arbitrary (nontrivial) abelian groups. Let $\Gamma$ be a nontrivial abelian group with identity $0$. A \textit{$\Gamma$-zonal labeling} of a plane graph $G$ is a labeling $\ell:V(G)\rightarrow \Gamma\setminus \{0\}$ such that for every region $f\in F(G)$ with boundary $B_f,$ $\sum_{v\in V(B_f)}\ell(v)=0\in \Gamma.$ We call this sum the \textit{induced label of $f$,} denoted $\ell(f).$

  Two additional conditions on $\Gamma$-zonal labelings were introduced in \cite{Bowling1}. We say a $\Gamma$-zonal labeling is \textit{generative} if the labels used generate $\Gamma.$ Such a restriction can be useful to consider. For example, the graph $K_4$ admits a $\Gamma$-zonal labeling for any group $\Gamma$ containing an element of order $3$; however, the only group for which it admits a generative $\Gamma$-zonal labeling (up to isomorphism) is $\mathbb{Z}_3$ (this follows from results in \cite{Bowling1} on wheel graphs).

 Our second restriction is the \textit{strong} condition. Let $G$ be a $2$-connected plane graph. A $\Gamma$-zonal labeling $\ell$ of $G$ is said to be strong if for all $v\in V(G),$ the degree of $v$ is equal to the additive order of the label $\ell(v)\in \Gamma.$ Symbolically, for all $v\in V(G), \deg(v)=\ord(\ell(v)).$ This not only draws a closer connection between the graph $G$ and the group $\Gamma$, but it also generalizes the key feature that connects zonal labelings to the Four Color Theorem in the $\mathbb{Z}_3$ case (for a more precise description of how, see \cite{Bowling1}).
 
 With the definition of zonality broadened to general abelian groups $\Gamma$, there are two primary questions:
 
 \begin{enumerate}
     \item Given a group $\Gamma$, for which graphs $G$ does $G$ admit a $\Gamma$-zonal labeling?
     \item Given a graph $G$, for which abelian groups $\Gamma$ does $G$ admit a $\Gamma$-zonal labeling?
\end{enumerate}

These questions can be extended by adding the generative and strong restrictions. Here, we will focus primarily on this second question of determining which groups are compatible with a given graph $G$. First, we will examine different families of graphs $\mathcal{G},$ and for each family we will explore for which groups $\Gamma$ there exists a graph $G\in \mathcal{G}$ admitting a $\Gamma$-zonal labeling (under various restrictions). Second, we will examine a tool that, given a specific graph $G$, can efficiently obtain information about the groups $\Gamma$ for which $G$ admits a $\Gamma$-zonal labeling.

 \section{A useful observation on composing functions}
Before we begin, we will present some important results on composing labelings with group homomorphisms. The first result appears in \cite{Bowling1}.

\begin{proposition}
\label{quotient}
Let $G$ be a plane graph, let $\Gamma, \Gamma'$ be abelian groups, and let $\alpha:\Gamma\rightarrow \Gamma'$ be a group homomorphism. If $\ell:V(G)\rightarrow \Gamma \setminus \ker(\alpha)$ is a $\Gamma$-zonal labeling, then $\alpha(\ell)=\alpha\circ \ell$ is a $\Gamma'$-zonal labeling. 
\end{proposition}

When our homomorphism is injective (and thus a \textit{monomorphism}), this result can be strengthened slightly.

\begin{corollary}
    Let $G$ be a plane graph, let $\Gamma, \Gamma'$ be abelian groups, and let $\alpha:\Gamma\rightarrow \Gamma'$ be a monomorphism. If $\ell:V(G)\rightarrow \Gamma\setminus \{0\}$ is a $\Gamma$-zonal labeling, then $\alpha(\ell)=\alpha\circ \ell$ is a $\Gamma'$-zonal labeling. 
\end{corollary}

In fact, in the case of strong $\Gamma$-zonal labelings, strongness is also preserved under monomorphisms, since these maps preserve the additive orders of elements. Therefore, we have the following corollary as well.
\begin{corollary}
    Let $G$ be a plane graph, let $\Gamma, \Gamma'$ be abelian groups, and let $\alpha:\Gamma\rightarrow \Gamma'$ be a monomorphism. If $\ell:V(G)\rightarrow \Gamma\setminus \{0\}$ is a strong $\Gamma$-zonal labeling, then $\alpha(\ell)=\alpha\circ \ell$ is a strong $\Gamma'$-zonal labeling. 
\end{corollary}

Lastly, it is useful to note that when our homomorphism is surjective, the generative property is preserved as well.

\begin{corollary}
    Let $G$ be a plane graph, let $\Gamma, \Gamma'$ be abelian groups, and let $\alpha:\Gamma\rightarrow \Gamma'$ be an epimorphism. If $\ell:V(G)\rightarrow \Gamma\setminus \ker(\alpha)$ is a generative $\Gamma$-zonal labeling, then $\alpha(\ell)=\alpha\circ \ell$ is a generative $\Gamma'$-zonal labeling. 
\end{corollary}
\begin{proof}
    Note that $\alpha(\ell)$ is $\Gamma$-zonal. Let $g'\in \Gamma'$. Then there is some $g\in \Gamma$ such that $\alpha(g)=g'.$ Furthermore, since $\ell$ is generative, $g$ is given by some integer linear combination $g=\sum_{v\in V(G)}k_v\ell(v)$. Thus, $g'=\alpha(g)=\sum_{v\in V(G)}k_v(\alpha\circ \ell)(v)$, and $\alpha(\ell)$ is generative.
\end{proof}

These basic observations will prove very useful in addressing both of our major questions, as we may often use a labeling in one group to derive a labeling in another.

\section{$\Gamma$-zonal labelings in 2-connected plane graphs}
We will begin with our first question: Given a family $\mathcal{G}$ of graphs, for which $\Gamma$ does there exist a plane graph $G_\Gamma\in \mathcal{G}$ having $\Gamma$-zonal labelings with various characteristics? Our attention will be focused on $2$-connected plane graphs. Specifically, we will examine the family $\mathcal{G}_{2,m}$ of $2$-connected plane multigraphs, and the family $\mathcal{G}_{2,s}$ of $2$-connected simple plane graphs. 

Our first object of study will be generative labelings, as there is a very general result that applies to both families. Note that if $\ell$ is a generative $\Gamma$-zonal labeling, then $\Gamma$ must be finitely generated, as otherwise the finite number of labels used in the labeling could not possibly generate $\Gamma$.

Observe that for $n\ge 3,$ the cycle $C_n$ belongs to $\mathcal{G}_{2,s}$, and therefore also to $\mathcal{G}_{2,m}$. The following was shown in \cite{Bowling1}.

\begin{theorem}
    Let $C_n$ be a cycle for $n\ge 3,$ and let $\Gamma \not\simeq\mathbb{Z}_2$ be a nontrivial abelian group. If $\Gamma$ is generated by at most $n-1$ elements, then $C_n$ is generative $\Gamma$-zonal.
\end{theorem}

Moreover, the even cycles $C_{2k}$ can be shown to be generative $\mathbb{Z}_2$-zonal by labeling every vertex $1$, as then the induced label of each region is $0$. We therefore have the following result.

\begin{theorem}
    Let $\Gamma$ be a nontrivial finitely generated abelian group. Then there exist infinitely many $2$-connected simple plane graphs that are generative $\Gamma$-zonal.
\end{theorem}
\begin{proof}
    First, let $\Gamma\simeq \mathbb{Z}_2.$ Then, the cycle $C_{2k}$ is generative $\mathbb{Z}_2$-zonal for all $k\ge 2.$

    Next, let $\Gamma\not\simeq \mathbb{Z}_2,$ such that $\Gamma$ is generated by a set of $m$ (non-identity) elements. Then the cycle $C_{n}$ is generative $\Gamma$-zonal for all $n>m.$

    Thus, in both cases, there are infinitely many cycles that are generative $\Gamma$-zonal, and therefore an infinite number of $2$-connected simple plane graphs that are generative $\Gamma$-zonal.
\end{proof}

We can use this result to generalize to $\Gamma$-zonal labelings as a whole.

\begin{corollary}
    Let $\Gamma$ be a nontrivial abelian group. Then there exist infinitely many $2$-connected simple plane graphs that are $\Gamma$-zonal.
\end{corollary}
\begin{proof}
    Choose an element $g\in \Gamma\setminus \{0\}.$ Then the group $\langle g \rangle$ is finitely generated, and thus there are infinitely many simple plane graphs that admit a generative $\langle g \rangle$-zonal labeling. Given such a graph $G_\Gamma$ and the generative $\langle g \rangle$-zonal labeling $\ell:V(G_\Gamma)\rightarrow \langle g \rangle\setminus \{0\},$ we compose with the natural monomorphism $\alpha:\langle g \rangle \rightarrow \Gamma$ to obtain a $\Gamma$-zonal labeling $\alpha(\ell)$ of $G_\Gamma$. Thus, each of the infinitely many graphs that are generative $\langle g \rangle$-zonal are also $\Gamma$-zonal. 
\end{proof}

Now that we have settled the question of the existence of generative $\Gamma$-zonal labelings in $2$-connected plane graphs, we now turn to generative strong $\Gamma$-zonal labelings in $2$-connected plane multigraphs. We will obtain a result that is equally general for this family of graphs. 

\begin{theorem}
\label{2connmulti}
    Let $\Gamma$ be a nontrivial abelian group. Then the following are equivalent.
    \begin{enumerate}
        \item There exists a $2$-connected plane multigraph that is generative strong $\Gamma$-zonal.
        \item There exist infinitely many $2$-connected plane multigraphs that are generative strong $\Gamma$-zonal.
        \item $\Gamma$ is finite.
    \end{enumerate}
\end{theorem}
\begin{proof}
    $3\implies 2)$ Let $\Gamma$ be a finite abelian group. Let $g_1, g_2, ...,g_n$ be a multiset of (nonzero) elements that generate $\Gamma$, and begin with the $2n$-cycle $$(u_1, v_1, u_2, v_2, ..., u_n, v_n, u_1).$$ (If there is only one generator $g_1,$ then count this generator twice so that we have at least a $4$-cycle.) Now, for each generator $g_i$, add $\ord(g_i)-2$ additional edges between $u_i$ and $v_i$, so that $\deg(u_i)=\deg(v_i)=\ord(g_i).$ Denote the resulting graph by $G_\Gamma$. To obtain a generative strong $\Gamma$-zonal labeling of $G_\Gamma,$ consider the labeling $\ell:V(G_\Gamma)\rightarrow \Gamma$ given by $\ell(u_i)=g_i$ and $\ell(v_i)=-g_i,$ for all $1\le i \le n.$ Now, each region $f$ of $G_\Gamma$ is either bounded by a $2$-cycle $(u_i, v_i, u_i),$ or is bounded by a cycle $(u_1, v_1, ..., u_n, v_n, u_1).$ In the first case, $\ell(f)=\ell(u_i)+\ell(v_i)=g_i-g_i=0.$ In the second case, $\ell(f)=\sum_{i=1}^n(\ell(u_i)+\ell(v_i))=\sum_{i=1}^n(g_i-g_i)=0.$ Thus, $\ell$ is a $\Gamma$-zonal labeling. Since $\deg(u_i)=\ord(g_i)=\ord(\ell(u_i)),$ and similarly for $v_i$, we see that $\ell$ is also strong. Lastly, since the labels used contain a set of generators for $\Gamma,$ the labeling is also generative.

     \begin{figure}[ht]\centering
\includegraphics[width=4in]{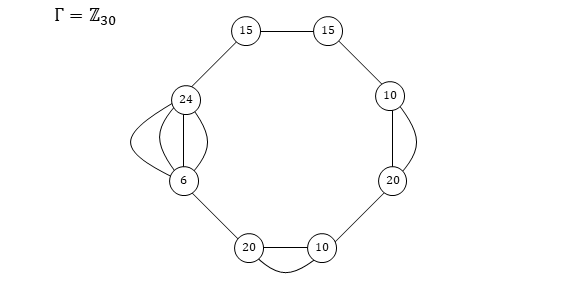}
\caption{A generative strong $\mathbb{Z}_{30}$-zonal labeling of a $2$-connected plane multigraph. This was formed using the multiset of generators $\{15, 10, 10, 6\}$} \label{z30Ex}
\end{figure}

    Note that the order of the graph can be made arbitrarily large by counting generators with higher multiplicity. Thus, there are infinitely many $2$-connected plane multigraphs that are generative strong $\Gamma$-zonal.

    $2\implies 1)$ This implication is immediate.

    $1\implies 3)$ Suppose that there exists a graph $G$ that admits a generative strong $\Gamma$-zonal labeling $\ell$ for some abelian group $\Gamma$. Then, $\Gamma$ is generated by the labels $\ell(V(G))$, each of which has finite order, and of which there are finitely many. Thus $\Gamma$ is finite.
\end{proof}

An example of a graph constructed using Theorem \ref{2connmulti} is given in Figure \ref{z30Ex}. Using this result, we now obtain a result on strong $\Gamma$-zonal labelings.

\begin{corollary}
    Let $\Gamma$ be an abelian group. Then the following are equivalent:
    \begin{enumerate}
        \item There exists a $2$-connected plane multigraph that is strong $\Gamma$-zonal.
        \item There exist infinitely many $2$-connected plane multigraphs that are strong $\Gamma$-zonal.
        \item $\Gamma$ contains an element $g$ with finite order $\ord(g)>1.$
    \end{enumerate} 
\end{corollary}
\begin{proof}
    $3\implies 2)$ Let $\Gamma$ be an abelian group with a nonzero element $g$ of finite order. Then the cyclic subgroup $\langle g\rangle$ is finite, and thus there are infinitely many $2$-connected plane multigraphs that admit a generative strong $\langle g \rangle$-zonal labeling. Given such a graph $G_\Gamma$ and the generative strong $\langle g \rangle$-zonal labeling $\ell:V(G_\Gamma)\rightarrow \langle g \rangle\setminus \{0\},$ we compose with the natural monomorphism $\alpha:\langle g \rangle \rightarrow \Gamma$ to obtain a strong $\Gamma$-zonal labeling $\alpha(\ell)$ of $G_\Gamma$. Thus, each of the infinitely many graphs that are generative strong $\langle g \rangle$-zonal is also strong $\Gamma$-zonal.   

    $2\implies 1)$ This is obvious.

    $1 \implies 3)$ Now, suppose that some $2$-connected plane multigraph $G$ admits a strong $\Gamma$-zonal labeling. Then, each label is a nonzero element with finite order, and therefore $\Gamma$ contains at least one nonzero element of finite order.
 \end{proof}

Next, we turn our attention more specifically to $\mathcal{G}_{2,s}$. In this case, there is an obstruction that prevents us from finding a graph $G_\Gamma\in \mathcal{G}_{2,s}$ admitting a strong $\Gamma$-zonal labeling altogether. We start here.

\begin{theorem}
\label{has235}
Let $\Gamma$ be an abelian group, and let $G$ be a $2$-connected simple plane graph. If $G$ admits a strong $\Gamma$-zonal labeling, then $\Gamma$ must have an element of order $2, 3,$ or $5$.   
\end{theorem}
\begin{proof}
    Assume by way of contradiction that $\Gamma$ has no such element, and that $G$ has a $\Gamma$-zonal labeling $\ell$. Note that $\Gamma$ also does not have any element of order $4$, as for any such element $g$, it follows that $\ord(2g)=2$. Note also that $G$ must have a vertex of degree $5$ or less. Since $G$ is $2$-connected, it also has no vertex of degree $1$. Therefore, $G$ has a vertex $v$ of degree $2, 3, 4,$ or $5$. Since $\ell$ is strong, then the additive order of $\ell(v)$ must be either $2, 3, 4,$ or $5$.  However, $\Gamma$ has no such element. This contradicts the existence of a strong $\Gamma$-zonal labeling, and therefore, $\Gamma$ must have an element of order $2, 3,$ or $5$.
\end{proof}

Using the graphs constructed in Theorem \ref{2connmulti} as a base, we can characterize the groups $\Gamma$ for which there is a graph in $\mathcal{G}_{2,s}$ with a generative strong $\Gamma$-zonal labeling.
\begin{theorem}
\label{2connsimple}
    Let $\Gamma$ be a finite abelian group. Then the following are equivalent:
    \begin{enumerate}
        \item There exists a $2$-connected simple plane graph that is generative strong $\Gamma$-zonal.
        \item There exist infinitely many $2$-connected simple plane graphs that are generative strong $\Gamma$-zonal.
        \item $\Gamma$ contains an element of order $2,3,$ or $5$.
    \end{enumerate}  
\end{theorem}

\begin{proof}
$1\implies 3)$ Suppose there is a $2$-connected plane simple graph that is strong $\Gamma$-zonal. By Theorem \ref{has235}, $\Gamma$ has an element of order $2, 3$, or $5$.

$3\implies 2)$
    Suppose that $\Gamma$ contains an element $g$ of order $2, 3,$ or $5$. By Theorem \ref{2connmulti}, there exists a $2$-connected plane multigraph $G_\Gamma$ with a generative $\Gamma$-zonal labeling $\ell$. We construct a new graph $G_\Gamma'$ that is simple and also has a generative strong $\Gamma$-zonal labeling.

    \begin{itemize}
        \item[]\textit{Case 1: $\ord(g)=2.$} Form $G_\Gamma'$ by subdividing each edge of $G$ twice, and form the labeling $\ell'$ by assigning each newly created vertex the label $g$ (and keeping all other labels the same). The resulting graph is simple, as no vertex originally in $G$ is adjacent to another vertex originally in $G$, and there are no parallel edges incident with newly created vertices. In addition, the induced label of each region in $G_\Gamma'$ is increased from the label of its corresponding region in $G_\Gamma$ by $2g=0$ for each edge originally on its boundary in $G_\Gamma$. Thus, the induced label of each region in $G_\Gamma'$ is the same as the corresponding region in $G_\Gamma$, which is $0$. Thus, $\ell'$ is zonal, and remains generative. Lastly, $\ell'$ is strong, as none of the original vertices changed degrees or labels, and each new vertex has degree $2$ and was assigned a label with additive order $2$. 

        \begin{figure}[ht]\centering
\includegraphics[width=4in]{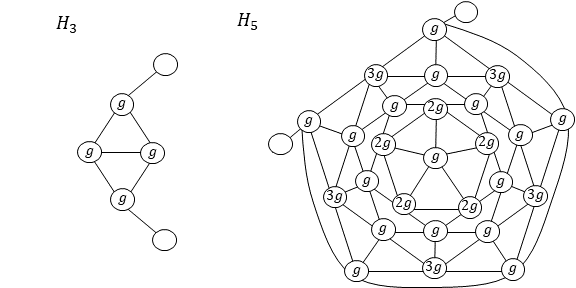}
\caption{Structures that replace edges for groups containing an element of order $3$ or $5$. The unlabeled vertices are the original endpoints of the edge.} \label{h3h5}
\end{figure}

        \item[]\textit{Case 2: $\ord(g)=3.$} Form $G_\Gamma'$ by replacing each edge of $G$ with the structure $H_3$ depicted in Figure \ref{h3h5}, and form the labeling $\ell'$ by assigning each newly created vertex the label $g$ (and keeping all other labels the same). The resulting graph is also simple. In addition, we can see that the induced label of each region in $G_\Gamma'$ (other than the regions interior to a copy of $H_3$) is increased from the label of its corresponding region in $G_\Gamma$ by $3g=0$ for each edge originally on its boundary in $G_\Gamma$ (as $G_\Gamma$ is $2-$connected and only one ``side" of $H_3$ is incident with a given region). In addition, the new regions formed as interior structures to $H_3$ also have an induced label of $0$. Thus, the induced label of each region in $G_\Gamma'$ is either the same as the corresponding region in $G_\Gamma$, or is $0$ as an interior region to a copy of $H_3$. Thus, $\ell'$ is zonal, and clearly still generative. Lastly, $\ell'$ is strong, as none of the original vertices changed degrees or labels, and each new vertex has degree $3$ and was assigned a label with additive order $3$. 

        \item[]\textit{Case 3: $\ord(g)=5.$} Form $G_\Gamma'$ by replacing each edge of $G$ with the structure $H_5$ depicted in Figure \ref{h3h5}, and form the labeling $\ell'$ by assigning each newly created vertex the label indicated in the figure (and keeping all other labels the same). Again, this is a simple graph. In addition, we can see that the induced label of each region in $G_\Gamma'$ (other than the regions interior to a copy of $H_5$) is increased from the label of its corresponding region in $G_\Gamma$ by either $g+3g+g=0$ or $g+g+g+g+g=0$ for each edge originally on its boundary in $G_\Gamma$. In addition, the new regions formed as interior structures to $H_5$ also have an induced label of $0$. Thus, the induced label of each region in $G_\Gamma'$ is either the same as the corresponding region in $G_\Gamma$, or is $0$ as an interior region to a copy of $H_5$. Thus, $\ell'$ is zonal, and is also generative and strong, as in the preceding cases.
    \end{itemize}

As there are infinitely many choices for the graph $G_\Gamma$, there are similarly infinitely many constructions $G_\Gamma'$. 

$2\implies 1)$ This is obvious.
\end{proof}

\begin{corollary}
    Let $\Gamma$ be an abelian group. Then the following are equivalent:
    \begin{enumerate}
        \item There exists a $2$-connected simple plane graph that is strong $\Gamma$-zonal.
        \item There exist infinitely many $2$-connected simple plane graphs that are strong $\Gamma$-zonal.
        \item $\Gamma$ contains an element of order $2,3,$ or $5$.
    \end{enumerate}  
\end{corollary}
\begin{proof}
    $3\implies 2)$ Let $g$ be an element of order $2, 3,$ or $5$. By Theorem \ref{2connsimple}, there exist infinitely many plane graphs that have a generative strong $\langle g \rangle$-zonal labeling. By composing with the natural monomorphism $\alpha:\langle g \rangle\rightarrow \Gamma$, any $2$-connected simple plane graph with a generative strong $\langle g \rangle$-zonal labeling also has a strong $\Gamma$-zonal labeling.

    $2\implies 1)$ This is obvious.

    $1\implies 3)$ This follows directly from Theorem \ref{has235}.
\end{proof}

We will now recast our existence results into compatibility claims between graphs in certain families and the groups for which a $\Gamma$-zonal labeling is feasible. Note that these claims are weaker than those proven above, but these formulations may be more useful in practice.

\begin{theorem}
\label{main1}
    The following are all true:
    \begin{enumerate}
        \item There exists a $2$-connected plane multigraph and a $2$-connected plane simple graph with a $\Gamma$-zonal labeling for any nontrivial abelian group $\Gamma$.
        \item Let $G$ be a $2$-connected plane multigraph or a $2$-connected plane simple graph with a generative $\Gamma$-zonal labeling. Then $\Gamma$ must be finitely generated.
        \item Let $G$ be a $2$-connected plane multigraph with a strong $\Gamma$-zonal labeling. Then $\Gamma$ must contain an element $g$ with finite order $|g|>1.$
        \item Let $G$ be a $2$-connected simple plane graph with a strong $\Gamma$-zonal labeling. Then $\Gamma$ contains an element of order $2, 3,$ or $5$.
        \item Let $G$ be a $2$-connected plane multigraph with a generative strong $\Gamma$-zonal labeling. Then $\Gamma$ must be finite.
        \item Let $G$ be a $2$-connected simple plane graph with a generative strong $\Gamma$-zonal labeling. Then $\Gamma$ must be finite and contain an element of order $2,3,$ or $5$.
    \end{enumerate}
\end{theorem}
Moreover, no stronger restrictions can be made than those listed above, as for each of the viable groups $\Gamma$ there exist infinitely many graphs $G$ in the given family that have the respective labeling.

\section{An algebraic perspective}

We will now move from studying families of graphs to studying individual graphs. While attempting to completely define for which $\Gamma$ a graph $G$ admits a $\Gamma$-zonal labeling would be challenging, we find that we can extract substantial information from routine algebraic calculations. 

Let $G$ be a plane graph with vertices $v_1,..., v_n$ and faces $f_1,...,f_r$. The \textit{face-vertex incidence matrix} of a plane graph $G$ is the matrix $M\in \mathbb{Z}^{r\times n}$ such that the entry $M_{ij}$ is $1$ if $v_j$ is on the boundary of $f_i,$ and $0$ otherwise. We will refer to this simply as a \textit{face incidence matrix} and let the vertex condition be understood. The face incidence matrix can be thought of as representing a $\mathbb{Z}$-linear transformation between the free $\mathbb{Z}$-modules $\mathbb{Z}^n$ and $\mathbb{Z}^r$. This matrix has appeared in various fields of mathematics, inluding the settings studied in \cite{Kaibel, Yan}.

We analyze this matrix using its Smith normal form (SNF). A matrix is \textit{unimodular} if its determinant is either $1$ or $-1$. These matrices are of interest because they are invertible over the ring $\mathbb{Z}$, and in fact over any ring $\mathbb{Z}_k$ where $k\ge 2.$ Given a matrix $M\in \mathbb{Z}^{r\times n}$, there exist unimodular matrices $U\in \gl_r(\mathbb{Z}),V\in \gl_n(\mathbb{Z})$ and a diagonal matrix $D\in \mathbb{Z}^{r\times n}$ such that $UMV=D.$ Furthermore, the diagonal entries $d_1, d_2,..., d_{\min(r,n)}$ of $D$ satisfy $d_1|d_2|...|d_{\min(r,n)}$. This equation $UMV=D$ is called the \textit{Smith normal form} decomposition. Several authors have studied labelings via analysis of the relevant incidence matrix \cite{Dickenstein, Jeurissen}. In addition, the Smith normal form has been used to analyze these incidence matrices \cite{Dickenstein, Grossman}. Here, we will see that the Smith normal form  of the face incidence matrix is not only useful in analyzing labelings, but in fact can be used to directly construct a zonal labeling when one exists. 

Given that our face incidence matrix represents a linear transformation, we will often be thinking of our groups as $R$-modules, where $R=\mathbb{Z}$ or $R=\mathbb{Z}_k$. A group $\Gamma$ is said to be \textit{annihilated} by $k$ if $kg=0$ for all $g\in \Gamma$. Note that any group $\Gamma$ is naturally a $\mathbb{Z}$-module, and if $k$ annihilates $\Gamma,$ the structure of $\Gamma$ as a group and as a $\mathbb{Z}_k$-module is identical. We will use $A_R$ to denote that the matrix $A$ is to be taken over the ring $R$. For additional module-theoretic terminology, see \cite{DummitFoote}. We begin by recasting zonal labelings into an entirely algebraic setting, determined by the face incidence matrix.

\begin{proposition}
\label{kernelProp}
    Let $G$ be a plane graph with face incidence matrix $M$, let $\Gamma$ be an abelian group, and let $\ell:V(G)\rightarrow \Gamma$ be a labeling of the vertices of $G$. Let the ring $R$ be either $\mathbb{Z}$ or $\mathbb{Z}_k$ for some $k$ that annihilates $\Gamma$. Let the standard basis of $R^r$ be given by $\tilde f_1,..., \tilde f_r$ corresponding to the faces $f_1, ..., f_r$ of $G$, let the standard basis of $R^n$ be given as $\tilde v_1,..., \tilde v_n$ corresponding to the vertices $v_1,..., v_n$ of $G$, and define $\tilde \ell:R^n\rightarrow \Gamma$ by $\tilde \ell(\tilde v_i)=\ell(v_i).$ Then $\ell:V(G)\rightarrow \Gamma$ is a $\Gamma$-zonal labeling if and only if $\im(M_R^T)\le \ker(\tilde \ell),$ and $\{\tilde v_1,..., \tilde v_n\}\cap \ker(\tilde \ell)=\emptyset$. 
\end{proposition}
\begin{proof}
    Begin by noting that the $j$th column of $M^T$ has ones in the rows corresponding to vertices on the boundary of $f_j$ and zeros elsewhere. Thus, $M_R^T(\tilde f_j)=\sum_{v\in B_{f_j}}\tilde v$. 

    Suppose that $\ell$ is a $\Gamma$-zonal labeling. Then $(\tilde \ell\circ M_R^T)( \tilde f_j)=\tilde \ell\left(\sum{_{v\in B_{\tilde f_j}}}\tilde v\right)=\sum{_{v\in B_{f_j}}}\ell(v)=0.$ Since the elements $\tilde f_1,..., \tilde f_r$ form a basis for $R^r$, this means that $\tilde \ell\circ M_R^T$ maps identically to $0$, and thus $\im(M_R^T)\le \ker(\ell').$ In addition, since $\ell$ is zonal, no vertex has a label of $0$, and thus $\tilde \ell(\tilde v_i)\ne 0$ for all $\tilde v_1,..., \tilde v_n.$ Thus, $\{\tilde v_1,..., \tilde v_n\}\cap \ker(\tilde \ell)=\emptyset$.

    Now suppose that $\im(M_R^T)\le \ker(\tilde \ell),$ and $\{\tilde v_1,..., \tilde v_n\}\cap \ker(\tilde\ell)=\emptyset$. Since $\tilde \ell\circ M_R^T$ is the zero function, this implies that $\sum_{v\in B_{f_j}}\tilde \ell(\tilde v)=\sum_{v\in B_{f_j}}\ell(v)=0$ for all $f_1,...,f_r$, and thus the induced face labels are all $0$. In addition, since $\{v_1,..., v_n\}\cap \ker(\ell')=\emptyset,$ this implies that $\ell(v_i)\ne 0$ for any $v_1,..., v_n.$ Therefore, $\ell$ is a $\Gamma$-zonal labeling.
\end{proof}
Given this identification, we will primarily consider our $\Gamma$-zonal labelings as mapping from $R^n$ as opposed to from $V(G)$ for the remainder of this section, with the basis elements of $R^n$ given as $v_1,..., v_n$ and the basis of $R^r$ given as $f_1,..., f_r$. Given a plane graph $G$ with face incidence matrix $M$, an abelian group $\Gamma$, and a ring $R$ such that $R=\mathbb{Z}$ or $R=\mathbb{Z}_k$ for some $k$ annihilating $\Gamma$, we say that $\ell:R^n\rightarrow \Gamma$ is a $\Gamma$-zonal labeling if $\im(M_R^T)\le \ker(\ell),$ and $\{v_1,..., v_n\}\cap \ker( \ell)=\emptyset$. 

There are several key $\Gamma$-zonal labelings where our group $\Gamma$ is isomorphic to $R^n/H$ for some submodule $H\le R.$ In these cases, we will give some refined conditions for $\Gamma$-zonal labelings.

\begin{proposition}
\label{quotientProp}
    Let $G$ be a plane graph with face incidence matrix $M$, let $R$ be either the ring $\mathbb{Z}$ or $\mathbb{Z}_k$ for some $k$, let $H\le R^n$, and let $L\in R^{n\times n}.$ Then the rows of $L$ give a $R^n/H$-zonal labeling $\ell:R^n\rightarrow R^n/H$ defined by assigning $\ell(v_i)$ as the coset containing the $i$th row of $L$ if and only if no row of $L$ belongs to $H$, and every row of $M_RL$ does belong to $H.$
\end{proposition}
\begin{proof}
Begin by noting that $L^T$ represents a transformation from $R^n\rightarrow R^n.$ Let $\pi_H$ be the canonical mapping $\pi_H:R^n\rightarrow R^n/H.$ It follows that $\ell(v_i)=(\pi_H\circ L^T)(v_i).$
    
    Now, suppose that $\ell$ is $R^n/H$-zonal. Since the rows of $L$ determine the cosets of the labels of the vertices, this implies that no row belongs to the coset $0+H$, which in this case implies that no row belongs to $H$. Next, note that by the zonality of $\ell,$ $$\pi_H\circ L^T\circ M_R^T=\ell\circ M_R^T=0.$$

    Thus, $\im(L^T\circ M_R^T)\le \ker(\pi_H)=H.$ Therefore each column of $L^TM_R^T$ belongs to $H$, and every row of $(L^TM_R^T)^T=M_RL$ belongs to $H$.

    Now suppose that no rows of $L$ belong to $H$, and every row of $M_RL$ does belong to $H.$ Then, for $v_1,...,v_n$, $L^T(v_i)$ doesn't belong to $H$, $\ell(v_i)=(\pi_H\circ L^T)(v_i)$ is not the coset $0+H,$ and thus $\ker(\ell)\cap \{v_1,...,v_n\}=\emptyset.$ Furthermore, since every row of $M_RL$ belongs to $H$, every column of $L^TM_R^T$ is in $H$, and thus $(\pi_H\circ L^T\circ M_R^T)(f_j)=0+H$ for each basis element $f_1,..., f_r$ of $R^r$. Thus, $\pi_H\circ L^T\circ M_R^T=\ell\circ M_R^T=0$ and $\im(M_R^T)\le \ker(\ell).$ We conclude that $\ell$ is $R^n/H$-zonal.
\end{proof}

Given a matrix $L$, we will say that the $R^n/H$-zonal labeling $\ell:R^n\rightarrow R^n/H$ that assigns to each vertex $v_i$ the coset containing the $i$th row of $L$ is \textit{represented by $L$.} 

We provide two further observations on our algebraic reframing.

\begin{proposition}
    Let $\ell:R^n\rightarrow \Gamma$ be a $\Gamma$-zonal labeling, where $\Gamma$ is a finitely generated abelian group and $R$ is either $\mathbb{Z}$ or $\mathbb{Z}_k$ for some $k$ which annihilates $\Gamma.$ Then $\ell$ is generative if and only if $\ell$ is surjective.
\end{proposition}
\begin{proof}
    Suppose that $\ell$ is generative. Then $\im(\ell)$ contains a set of generators of $\Gamma$, and $\ell$ is surjective.

    Now suppose that $\ell$ is surjective, and let $g\in \Gamma$ be any element. Then there is some $v\in R^n$ such that $\ell(v)=g.$ Since $v=r_1v_1+r_2v_2+...+r_nv_n$ for scalars $r_1,...,r_n$, it follows by $R$-linearity that $g=r_1\ell(v_1)+r_2\ell(v_2)+...+r_n\ell(v_n)$, and $g$ is generated by the labels used in $\ell.$ Thus, $\ell$ is generative.
\end{proof}

\begin{proposition}
    Let $G$ be a plane graph with face incidence matrix $M$, let $R$ be either the ring $\mathbb{Z}$ or $\mathbb{Z}_k$ for some $k$, let $H\le R^n$, let $L\in R^{n\times n}$ represent a $R^n/H$-zonal labeling $\ell:R^n\rightarrow R^n/H,$ and let $K\in R^{n\times n}$ represent a linear transformation $R^n/H\rightarrow R^n/KH.$ Then the composition $K\circ \ell$ mapping from $R^n\rightarrow R^n/KH$ is given by the rows of $LK^T$.
\end{proposition}
\begin{proof}
    Note that the $i$th column of $KL^T$ is exactly the image of the $i$th column of $L^T$ under the transformation $K$. Taking the transpose, the $i$th row of $LK^T$ is the image of the $i$th row of $L$ under the transformation $K$. Since the $i$th row of $L$ is an element belonging to the coset $\ell(v_i),$ it follows that the $i$th row of $LK^T$ is an element belonging to the coset of $(K\circ \ell)(v_i).$ Therefore, the composition $K\circ \ell$ is given by the rows of $LK^T$. 

\end{proof}

We now use this algebraic framework to study zonality through the Smith normal form. 

\section{An algebraic tool for studying zonality}
Note that given $M_R^T:R^r\rightarrow R^n$, we have the exact sequence
$$0\rightarrow \ker(M_R^T)\rightarrow R^r\xrightarrow{M_R^T} R^n\rightarrow \operatorname{coker}(M_R^T)\rightarrow 0.$$

From this, we deduce that $\operatorname{coker}(M^T_R)=R^n/\operatorname{im}(M_R^T)$. This is isomorphic to $R^n/\im(D_R^T)$, which is in turn isomorphic to $\bigoplus_{i=1}^nR/\left<d_i\right>$, where $d_i$ is the $i$th diagonal element in $D_R^T$ from the SNF decomposition (or, for $r< i\le n, d_i=0)$. 

We proceed with our following theorem, which allows us to determine whether $G$ admits a $\Gamma$-zonal labeling for some $\Gamma$, by computing the Smith normal form decomposition of $M.$

\begin{theorem}
\label{main2}
    Let $G$ be a plane graph, and let $M$ be the face incidence matrix with Smith normal form decomposition $UMV=D.$ Let $R$ be either $\mathbb{Z}$ or $\mathbb{Z}_k$. Let $\pi:R^n\rightarrow \coker(M_R^T)$ map each basis element of $R^n$ to the corresponding element $e_i\in \coker(M_R^T),$ and let $\varphi:R^n\rightarrow R^n/\im(D_R^T)$ be represented by the rows of $V_R.$ The following are equivalent:
    \begin{enumerate}
        \item $G$ admits a $\Gamma$-zonal labeling for some abelian $\Gamma$. (In the case where $R=\mathbb{Z}_k$, $\Gamma$ must be a group that is annihilated by $k$.)
        \item The rows of $M_R$ do not generate a standard basis element of $R^n$.
        \item No row of $V_R$ is in $\im(D_R^T).$
        \item $\pi$ is a generative $\coker(M_R^T)$-zonal labeling.
        \item $\varphi$ is a generative $R^n/\im(D_R^T)$-zonal labeling.
    \end{enumerate}

    Now let $\Gamma$ be abelian, and let $\ell:R^n\rightarrow \Gamma$ be a particular $\Gamma$-zonal labeling (where, if $R=\mathbb{Z}_k,$ then $k$ annihilates $\Gamma$). Then,
    \begin{enumerate}[label=(\alph*)]
        \item $\ell=\ell'\circ \pi$ for a homomorphism $\ell':\coker(M_R^T)\rightarrow \Gamma.$
        \item $\ell=\ell''\circ \varphi$ for a homomorphism $\ell'':R^n/\im(D_R^T)\rightarrow \Gamma$.
    \end{enumerate}
\end{theorem}

\begin{proof}
    $1\implies 2$ Let $G$ admit a $\Gamma$-zonal labeling $\ell:R^n\rightarrow \Gamma$ for some $\Gamma$ (with $k$ annihilating $\Gamma$ if $R=\mathbb{Z}_k$). By Proposition \ref{kernelProp}, $\im(M^T_R)\le \ker(\ell)$ and $\{v_1,..., v_n\}\cap \ker(\ell)=\emptyset.$ Thus, $\{v_1,..., v_n\}\cap \im(M_R^T)=\emptyset,$ and the columns of $M_R^T$ (equivalently, the rows of $M_R$) do not generate a standard basis element of $R^n.$

    $2\implies 4$ Suppose that the rows of $M_R$ do not generate a standard basis element of $R^n.$ Consider the labeling $\pi:R^n\rightarrow \coker(M_R^T),$ noting that $\coker(M_R^T)=R^n/\im(M_R^T).$ Note that this labeling can be represented by $I_n.$ Since the columns of $M_R^T$ do not generate a standard basis element of $R^n$, and since the rows of $I_n$ are exactly the standard basis elements, no row of $I_n$ is in $\im(M_R^T).$ Furthermore, every row of $M_RI_n=M_R$ is clearly in $\im(M_R^T)$. Thus, by Proposition \ref{quotientProp}, $\pi$ is a $\coker(M_R^T)$-zonal labeling. Clearly $\pi$ is surjective and thus generative.

    $4\implies 5$ Suppose that $\pi$ is a generative $\coker(M_R^T)$-zonal labeling. Note that $V_R^TM_R^TU_R^T=D_R^T$. Since $V,U$ are unimodular, $V_R, U_R$, and their transposes are all invertible, and thus $M_R^TU_R^T=V^{-T}_RD_R^T$. Again, by the unimodularity of $U_R^T,$ $\im(M_R^T)=\im(V_R^{-T}D_R^T).$ Therefore, we may use the invertible map $V_R^{T}$ to map from $R^n/\im(M_R^T)$ to $R^n/\im(D_R^T).$ Now consider the composition $V_R^T\circ \pi:R^n\rightarrow R^n/\im(D_R^T).$ Since $\pi$ is $\coker(M_R^T)$-zonal and $V_R^T$ is an isomorphism, it follows that $V_R^T\circ \pi$ is $R^n/\im(D_R^T)$-zonal. Noting that the map $V_R^T\circ \pi$ is given by $I_nV_R=V_R$, this map $V_R^T\circ \pi$ is in fact $\varphi.$ Since both $\pi$ and $V_R^T$ are surjective, so is $\varphi$, and thus $\varphi$ is generative. 

    $5\implies 3$ Suppose that $\varphi$ is a generative $R^n/\im(D_R^T)$-zonal labeling. Then by Proposition \ref{quotientProp}, no row of $V_R$ belongs to $\im(D_R^T).$ 

    $3\implies 5$ Suppose that no row of $V_R$ belongs to $\im(D_R^T),$ and consider the labeling $\varphi:R^n\rightarrow R^n/\im(D_R^T)$ given by the rows of $V_R.$ Note that $V_R^TM_R^T=D_R^TU_R^{-T}.$ Thus, each column of $V_R^TM_R^T$ is in the image of $V_R^TM_R^T=D_R^TU_R^{-T},$ and thus in the image of $D_R^T.$ It follows that each row of $(V_R^TM_R^T)^T=M_RV_R$ is in the image of $D_R^T.$ By Proposition \ref{quotientProp}, $\varphi$ is zonal. Since $V_R$ is unimodular, $V_R$ and $V_R^T$ are surjective (in fact, isomorphisms), and thus the rows of $V_R$ generate $R^n$. Thus, $\varphi$ is surjective and thus a generative labeling.

    $5\implies 1$ $\varphi:R^n\rightarrow R^n/\im(D_R^T)$ is the desired $\Gamma$-zonal labeling. \\

    Now suppose that $\ell:R^n \rightarrow \Gamma$ is a particular $\Gamma$-zonal labeling.

    \begin{enumerate}[label=(\alph*)]
    \item Since $\ell$ is a $\Gamma$-zonal labeling, by Proposition \ref{kernelProp}, $\im(M_R^T)\le \ker(\ell).$ Thus, since $\pi$ is the canonical map from $R^n\rightarrow R^n/\im(M_R^T),$ there is a homomorphism $\ell':R^n/\im(M_R^T)\rightarrow \Gamma$ such that $\ell=\ell'\circ \pi.$  
    \item Let $\ell''=\ell'\circ V_R^{-T}.$ Then $\ell''\circ \varphi=\ell'\circ V_R^{-T}\circ V_R^T\circ \pi=\ell'\circ \pi=\ell,$ and we have that $\ell'':R^n/\im(D_R^T)\rightarrow \Gamma$ is a homomorphism such that $\ell=\ell''\circ \varphi.$
    \end{enumerate}

\end{proof}
We summarize the previous theorem in the following commutative diagram.

\[\begin{tikzcd}
	&&&& {R^n/\operatorname{im}(D_R^T)} & {\bigoplus_{i=1}^nR/\langle d_i\rangle} \\
	&&&& {\operatorname{coker}(M_R^T)} & {R^n/\operatorname{im}(M_R^T)} \\
	{R^r} && {R^n} && \Gamma \\
	&&&& 0
	\arrow["\simeq", tail reversed, from=1-5, to=1-6]
	\arrow["{V_R^{-T}}"', curve={height=12pt}, from=1-5, to=2-5]
	\arrow["{V_R^T}"', curve={height=12pt}, from=2-5, to=1-5]
	\arrow["{=}", tail reversed, from=2-5, to=2-6]
	\arrow["{\ell'}", from=2-5, to=3-5]
	\arrow["{M_R^T}", from=3-1, to=3-3]
	\arrow[curve={height=12pt}, from=3-1, to=4-5]
	\arrow["\varphi", curve={height=-30pt}, two heads, from=3-3, to=1-5]
	\arrow["\pi", two heads, from=3-3, to=2-5]
	\arrow["\ell", from=3-3, to=3-5]
	\arrow[hook, from=4-5, to=3-5]
\end{tikzcd}\]

In particular, if $\ell$ is zonal, then the rows of $M_R$ do not generate a basis element of $R^n.$ Thus, the map $\pi$ and the composition $\varphi$ are also zonal. Because $\ell$ is zonal, $\im(M_R^T)\le \ker(\ell),$ and we obtain the induced map $\ell'$ in the diagram. 

\begin{corollary}
\label{main3}
    Let $G$ be a plane graph with face incidence matrix $M$, let $\Gamma$ be an abelian group, and let $R$ be either $\mathbb{Z}$ or $\mathbb{Z}_k$ for some $k$ which annihilates $\Gamma$. If there exists a generative $\Gamma$-zonal labeling $\ell:R^n\rightarrow \Gamma$, then $\Gamma$ is isomorphic to a quotient group of $\coker(M_R^T),$ (or equivalently, isomorphic to a quotient group of $R^n/\im(D_R^T)$).
\end{corollary}
\begin{proof}
    Let $\ell$ be a generative $\Gamma$-zonal labeling. Then $\ell$ is surjective. Thus, since $\ell=\ell'\circ \pi$, $\ell':\coker(M^T)\rightarrow \Gamma$ is surjective, and $\Gamma\simeq \coker(M_R^T)/\ker(\ell').$ A similar argument using $\ell=\ell''\circ \varphi$ implies that $\Gamma\simeq (R^n/\im(D_R^T))/\ker(\ell'').$
\end{proof}

Note that it is entirely possible that $\Gamma$ is isomorphic to a quotient group of $\coker(M^T)$ and yet there is no $\Gamma$-zonal labeling (generative or otherwise). We will see an example of this case in the next section.

Theorem \ref{main2} and Corollary \ref{main3} have the following consequences:
\begin{enumerate}
    \item We can determine whether there exists a $\Gamma$ such that $G$ admits a $\Gamma$-zonal labeling simply by checking a divisibility rule along the rows of $V$. Indeed, for each row, we check whether the $i$th entry is divisible by $d_i$. If this is the case for every entry of the row, then that row is in the image of $D^T$, and there does not exist any $\Gamma$-zonal labeling. Otherwise, the rows of $V$ immediately give the zonal labeling $\varphi$.
    \item  We can determine whether there exists a $\Gamma$ such that $G$ admits a $\Gamma$-zonal labeling and such that $k$ annihilates $\Gamma$ by checking a modified divisibility rule. Let $R=\mathbb{Z}_k$. For each row of $V$, we check whether $\gcd(d_i,k)$ divides the $i$th entry. If this is the case for every entry, then the corresponding row of $V_R$ is in $\im(D_R^T),$ and there does not exist a $\Gamma$-zonal labeling for such $\Gamma$. Otherwise, we again have the labeling $\varphi.$ 
    \item Given a group $\Gamma$, if we want to find a $\Gamma$-zonal labeling, we can begin with $\varphi,$ and then any $\Gamma$-zonal labeling $\ell:R^n\rightarrow \Gamma$ will be given as the composition of $\varphi$ with some homomorphism $\ell':R^n/\im(D_R^T)\rightarrow \Gamma$. Furthermore, in forming a $\Gamma$-zonal labeling this way, we need only check that $\varphi(v_i)\notin \ker(\ell')$ for any vertex $v_i$, or equivalently, that $\ell'\circ \varphi$ does not assign the label $0$ to any vertex.
\end{enumerate}

\section{Application of the Smith normal form approach}
Here, we will use the tools developed in the previous section to examine some specific graphs. Although these graphs are primarily illustrative, they will illustrate several ways that the Smith normal form tools can be used.

We begin by examining the wheel $W_8$. In \cite{Bowling1}, it was shown that $W_8$ admits a generative $\Gamma$-zonal labeling if and only if $\Gamma$ can be generated by two nonzero elements $a,b\in \Gamma$ such that $a\ne -b$ and the additive order of $a+b$ divides $4$. Here, we will obtain similar results using the Smith normal form of the face incidence matrix.

 \begin{figure}[ht]\centering
\includegraphics[width=4in]{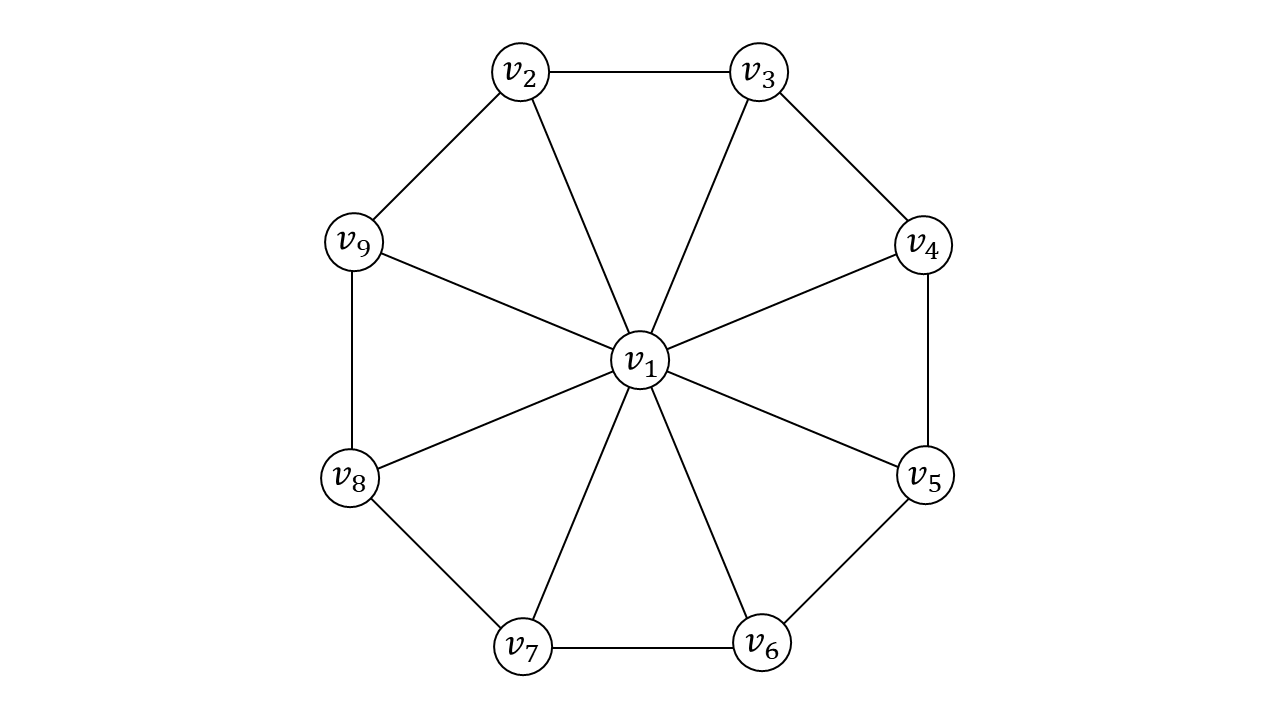}
\caption{The graph $W_8=C_8\vee K_1$, labeled to correspond with its face incidence matrix.} \label{w8Labeled}
\end{figure}

Consider the graph $W_8$ as drawn in Figure \ref{w8Labeled}. The face incidence matrix is given below, where the $i$th column corresponds to the vertex $v_i.$

\[
M=
\begin{bmatrix}
1&1&1&0&0&0&0&0&0\\
1&0&1&1&0&0&0&0&0\\
1&0&0&1&1&0&0&0&0\\
1&0&0&0&1&1&0&0&0\\
1&0&0&0&0&1&1&0&0\\
1&0&0&0&0&0&1&1&0\\
1&0&0&0&0&0&0&1&1\\
1&1&0&0&0&0&0&0&1\\
0&1&1&1&1&1&1&1&1
\end{bmatrix}
\]

When decomposed into Smith normal form as $D=UMV,$ the matrix $D$ has diagonals $d_i=1$ for $1\le i  \le 7$, $d_8=4,$ and $ d_9=0.$ The matrix $V$ is given below.

\[
V=
\begin{bmatrix}
0&0&0&0&0&0&0&-1&0\\
0&0&0&0&0&0&0&0&1\\
0&0&1&0&0&0&0&1&-1\\
0&0&-1&1&0&0&0&0&1\\
0&0&1&-1&1&0&0&1&-1\\
0&0&-1&1&-1&1&0&0&1\\
0&0&1&-1&1&-1&1&1&-1\\
0&1&-1&1&-1&1&-1&0&1\\
1&-1&1&-1&1&-1&1&1&-1
\end{bmatrix}
\]

To determine whether a zonal labeling exists for some target group, we determine whether any of the rows of $V$ are in the image of $D^T.$ This can be done by determining whether there is any row where $d_i$ divides the $i$th entry of that row for all entries. As $d_i=1$ for $1\le i\le 7,$ the only columns of significance are the last two. In this case, we see that no row has both an eighth entry divisible by $4$ and a final entry equal to $0$. Thus, a $\Gamma$-zonal labeling exists for some $\Gamma$. As $\mathbb{Z}^9/\im(D^T)\simeq \mathbb{Z}_4\oplus \mathbb{Z},$ we see that the generative $\mathbb{Z}_4\oplus \mathbb{Z}$-zonal labeling $\varphi:V(W_8)\rightarrow \mathbb{Z}_4\oplus \mathbb{Z}$ can be given by the last two columns of $V$, where the second to last column is taken to be in $\mathbb{Z}_4.$ The labeling is given as follows:

\[
\renewcommand{\arraystretch}{1.2}
\begin{array}{r@{{}={}}l@{\qquad}
                r@{{}={}}l@{\qquad}
                r@{{}={}}l}
\varphi(v_1) & (3,0)  &
\varphi(v_2) & (0,1)  &
\varphi(v_3) & (1,-1) \\
\varphi(v_4) & (0,1)  &
\varphi(v_5) & (1,-1) &
\varphi(v_6) & (0,1)  \\
\varphi(v_7) & (1,-1) &
\varphi(v_8) & (0,1)  &
\varphi(v_9) & (1,-1)
\end{array}
\]

We can draw several further conclusions. First, there is no generative $\mathbb{Z}_4\oplus \mathbb{Z}\oplus \mathbb{Z}$-zonal labeling of $W_8$, for instance, because this is not a quotient group of $\mathbb{Z}_4\oplus \mathbb{Z}$.

Additionally, there is no $\Gamma$-zonal labeling of $W_8$ where $\Gamma$ is annihilated by some odd $k\ge 3$. Recall that any $\Gamma$-zonal labeling $\ell:V(W_8)\rightarrow \Gamma$ is the composition of $\varphi$ with some homomorphism $\ell':\mathbb{Z}_4\oplus \mathbb{Z}\rightarrow \Gamma$. Since $\ord(\varphi(v_1))=4,$ the additive order of $\ell(v_1)$ divides $4$. If $\Gamma$ is annihilated by some odd $k$, then $\Gamma$ has no elements of even order, and thus $\ell(v_1)=\ell'(\varphi(v_1))=0,$ contradicting the zonality of $\ell.$ 

Furthermore, since $\ord(\ell(v_1))$ divides $4$ for any $\Gamma$-zonal labeling, there is no $\mathbb{Z}$-zonal labeling of $W_8$, as no nonzero element of $\mathbb{Z}$ has finite order. In addition, since $\deg(v_1)=8,$ we also conclude that there is no strong $\Gamma$-zonal labeling of $W_8.$ In fact, it was shown in \cite{Bowling1} that there is no strong $\Gamma$-zonal labeling of $W_n$ for any $n \ge 4$. 

On the other hand, there is a generative $\mathbb{Z}_2\oplus\mathbb{Z}_2$-zonal labeling $\varphi_2:V(W_8)\rightarrow \mathbb{Z}_2\oplus \mathbb{Z}_2$, which can be given by taking the last two columns of $V$ and reducing them modulo $2$. However, there is not a $\mathbb{Z}_2$-zonal labeling of $W_8.$ To see why, note that such a labeling $\ell:V(W_8)\rightarrow \Gamma$ would necessarily be the composition of $\varphi_2$ with some homomorphism $\ell':\mathbb{Z}_2\oplus \mathbb{Z}_2\rightarrow\mathbb{Z}_2$. There are exactly three such nonzero homomorphisms, given by $(a,b)\mapsto a, (a,b)\mapsto b$, and $(a,b)\mapsto a+b.$ Given that $(1,0), (0,1)$ and $(1,1)$ are all used as labels in $\varphi_2,$ at least one of these labels must belong to the kernel of $\ell',$ contradicting that $\ell=\ell'\circ \varphi$ is a $\mathbb{Z}_2$-zonal labeling. From a graph theory perspective, we can simply observe that no triangular region can have an induced label of $0$ if the only available label is $1\in \mathbb{Z}_2.$

\begin{figure}[ht]\centering
\includegraphics[width=4in]{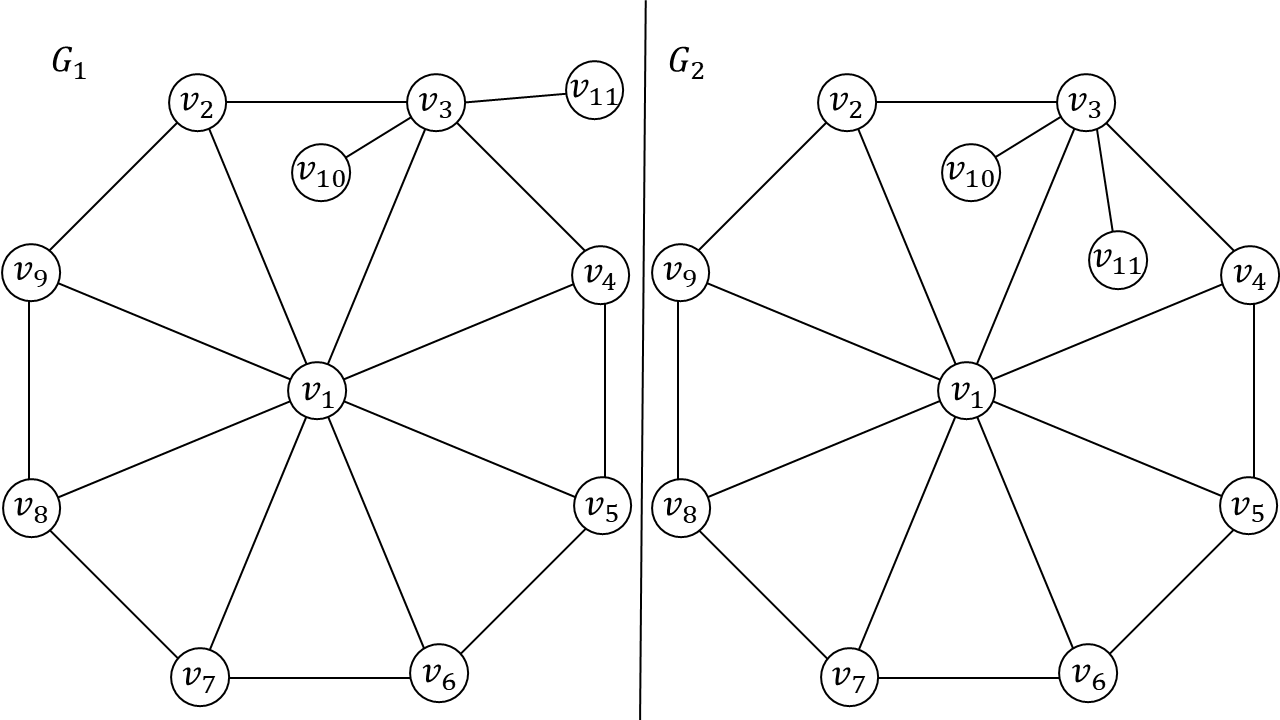}
\caption{Two embeddings of the graph formed by adding two pendant edges to $v_3.$} \label{g1g2SNFex}
\end{figure}

Next, consider a new plane graph $G$ formed by attaching two pendant edges to the vertex $v_3.$ We will examine the two embeddings $G_1, G_2$ presented in Figure \ref{g1g2SNFex}.

For brevity, the face incidence matrix of the graphs will be omitted here, although they can be obtained by appending columns to the face incidence matrix of $W_8$ corresponding to the vertices $v_{10}, v_{11}$. In each case, the diagonals of the matrix $D$ from Smith normal form are $d_i=1$ for $1\le i\le 9.$ As this matrix has more vertices than regions, we will append two extra entries $d_{10},d_{11}$ each equal to $0$ for convenience; this corresponds to each $d_i$ representing the element in $\mathbb{Z}^n$ containing $d_i$ in the $i$th position and zeros elsewhere. The agreement of the diagonal entries is coincidental; this is not true for all embeddings of $G$.

For the graph $G_1,$ the matrix $V$ from the Smith normal form is given below. 

\[
V=\begin{bmatrix}
0&0&0&0&0&0&0&0&0&1&0\\
0&0&0&0&0&0&0&0&0&0&1\\
0&0&-1&1&-1&1&-1&1&1&-1&-1\\
0&0&1&-1&1&-1&1&0&-1&0&1\\
0&0&-1&1&-1&1&0&0&1&-1&-1\\
0&0&1&-1&1&0&0&0&-1&0&1\\
0&0&-1&1&0&0&0&0&1&-1&-1\\
0&0&1&0&0&0&0&0&-1&0&1\\
0&0&0&0&0&0&0&0&1&-1&-1\\
0&1&1&-1&1&-1&1&-1&-1&0&0\\
1&0&0&-1&0&-1&0&-1&-1&4&0
\end{bmatrix}
\]

Examining the tenth row, the last two entries are $0$, indicating that this row belongs to $\im(D^T).$ Thus, there is no $\Gamma$-zonal labeling of $G_1.$ To justify this claim graph-theoretically, the triangular regions can be used to show that, in a $\Gamma$-zonal labeling, $\ell(v_4)=\ell(v_6)=\ell(v_8)=\ell(v_2).$ Thus, $\ell(v_3)+\ell(v_4)+\ell(v_1)=\ell(v_3)+\ell(v_2)+\ell(v_1)=0.$ Thus, for the region $f$ with $v_1, v_2, v_3, v_{10}$ on its boundary, we have that $\ell(f)=\ell(v_{10}),$ which is impossible in a $\Gamma$-zonal labeling as $\ell(f)=0$ and $\ell(v_{10})\ne 0$. 

Finally, we turn our attention to $G_2$. The matrix $V$ from Smith normal form is given below:

\[
V=\begin{bmatrix}
0&0&0&0&0&0&0&0&0&1&0\\
0&0&0&0&0&0&0&0&0&0&1\\
0&0&-1&0&-1&0&-1&0&1&3&-1\\
0&0&1&-1&1&-1&1&-1&0&0&1\\
0&0&-1&1&-1&1&0&1&0&-1&-1\\
0&0&1&-1&1&0&0&-1&0&0&1\\
0&0&-1&1&0&0&0&1&0&-1&-1\\
0&0&1&0&0&0&0&-1&0&0&1\\
0&0&0&0&0&0&0&1&0&-1&-1\\
0&1&1&0&1&0&1&0&-1&-4&0\\
1&0&0&1&0&1&0&1&-1&-4&0
\end{bmatrix}
\]

As there is no row where the last two entries are both $0$, no row belongs to $\im(D^T),$ and the rows of $V$ give rise to the $\mathbb{Z}\oplus \mathbb{Z}$-zonal labeling $\varphi:V(G_2)\rightarrow \mathbb{Z}\oplus \mathbb{Z}$:

\vspace{10pt}

\[
\renewcommand{\arraystretch}{1.2}
\begin{array}{l@{{}={}}l@{\qquad}
                l@{{}={}}l@{\qquad}
                r@{{}={}}l}
\varphi(v_1)     & (1,0)   &
\varphi(v_2)     & (0,1)   &
\varphi(v_3)     & (3,-1)  \\
\varphi(v_4)     & (0,1)   &
\varphi(v_5)     & (-1,-1) &
\varphi(v_6)     & (0,1)   \\
\varphi(v_7)     & (-1,-1) &
\varphi(v_8)     & (0,1)   &
\varphi(v_9)     & (-1,-1) \\
\varphi(v_{10})  & (-4,0)  &
\varphi(v_{11})  & (-4,0)  &
\multicolumn{2}{c}{}
\end{array}
\]

\vspace{10pt}

However, if one only considers groups $\Gamma$ that are annihilated by $2$ or $4,$ then there is no $\Gamma$-zonal labeling. To see why, note that when $R=\mathbb{Z}_k$, then $\gcd(d_{10},k)=\gcd(d_{11}, k)=k.$ In the tenth and eleventh rows of $V$, the last two entries are both divisible by $4$, and thus also by $2$. Therefore, these rows belong to $\im(D_R^T),$ and thus there is no $\Gamma$-zonal labeling where $\Gamma$ is annihilated by either $2$ or $4$.

In fact, there is a generative $\Gamma$-zonal labeling of $G_2$ for any other quotient group of $\mathbb{Z}\oplus \mathbb{Z}$. We divide into cases:
\begin{itemize}
    \item[] \textit{Case 1: $\Gamma\simeq \mathbb{Z}$}. Form $\ell_1=\ell_1'\circ \varphi,$ where $\ell_1':\mathbb{Z}\oplus \mathbb{Z}\rightarrow \Gamma$ is given by $\ell_1'(a,b)=a+b.$ 
    \item[] \textit{Case 2: $\Gamma\simeq \mathbb{Z}\oplus \mathbb{Z}_k$.} Form $\ell_2=\ell_2'\circ \varphi,$ where $\ell_2':\mathbb{Z}\oplus \mathbb{Z}\rightarrow \Gamma$ is given by $\ell_2'(a,b)=(a, b \pmod k).$ 
    \item[] \textit{Case 3: $\Gamma\simeq \mathbb{Z}_k,$ where $k\ge 3, k\ne4$.} Form $\ell_3=\ell_3'\circ \varphi,$ where $\ell_3':\mathbb{Z}\oplus \mathbb{Z}\rightarrow \Gamma$ is given by $\ell_3'(a,b)=a+b \pmod k.$ 
    \item[] \textit{Case 4: $\Gamma\simeq \mathbb{Z}_j\oplus \mathbb{Z}_k$, $j,k\ge 2,$ $j\not \in \{2,4\}.$} Form $\ell_4=\ell_4'\circ \varphi,$ where $\ell_4':\mathbb{Z}\oplus \mathbb{Z}\rightarrow \Gamma$ is given by $\ell_4'(a,b)=(a \pmod j, b \pmod k)$. 
\end{itemize}
In all cases, $\ell_i'$ is a surjection, and there is no vertex $v\in V(G_2)$ such that $\varphi(v)\in \ker(\ell_i').$ Thus, each of the compositions $\ell_i=\ell_i'\circ \varphi$ are generative $\Gamma$-zonal for the respective $\Gamma.$

\section{Concluding remarks}
First, consider Theorem \ref{main1}. The $2$-connected condition was chosen because strong zonal labelings are only defined for $2$-connected graphs, as this is the environment in which their motivating properties hold. However, one could ask about $3$-connected graphs instead. Alternatively, one could explore these questions from the dual perspective of cozonal labelings. Let $X_v$ denote the set of regions incident with $v$. A \textit{$\Gamma$-cozonal labeling} is a labeling $\ell:F(G)\rightarrow \Gamma\setminus \{0\}$ such that for all $v\in V(G), \sum_{f\in X_v}\ell(f)=0$. There are analogous definitions for generative and strong $\Gamma$-cozonal labelings, as described in \cite{Bowling1}.

The analogous cozonal questions are substantially more complex than their zonal counterparts, particularly when considering the generative condition. However, there are much deeper tools that can be developed from the cozonal perspective, which can then be used to prove strong claims about $3$-connected plane graphs as well. These topics will be explored in a future publication.

Now, let us consider the results of Theorem \ref{main2} and Corollary \ref{main3}. As we saw in our examples, this tool is especially well-suited for studying labelings of individual graphs. Given a graph $G$, the labeling $\varphi$ constructed using Smith normal form can be used to determine many necessary conditions for arbitrary $\Gamma$-zonal labelings; in particular, $\varphi$ often establishes restrictions on the additive orders of elements. These restrictions can be invaluable in determining the groups $\Gamma$ for which a graph admits a $\Gamma$-zonal labeling, as well as in determining whether a strong $\Gamma$-zonal labeling can be found. Moreover, to find various $\Gamma$-zonal labelings, one need not necessarily find the labeling via pure trial and error; instead, one can begin with $\varphi$ and then find some homomorphism $\ell'$ such that no label $\varphi(v_i)$ belongs to $\ker(\ell').$ Nevertheless, this too involves some degree of trial and error.

As for the analysis of families of graphs, it is unclear if the Smith normal form approach can work effectively in this context. In the case where a family of graphs has a standard enough face incidence matrix that a Smith normal form decomposition can be derived uniformly, this may be a very powerful tool. Even if this is not the case, examining the labeling $\varphi$ on specific instances can help understand the structure of a $\Gamma$-zonal labeling, or otherwise the obstacles that prevent one from being constructed.

\section{Tool and computational resource statement}\label{Disclosure}The authors support the Leiden Declaration on Artificial Intelligence and Mathematics. An informal conversation with ChatGPT provided inspiration for this paper. ChatGPT was also used to proofread the manuscript for grammar, spelling, and punctuation. We take full responsibility for the mathematical content of this paper.


\begin{thebibliography}{99}

\bibitem{Bowling}
Bowling, A., \textit{Zonality in Graphs}. Ph.D. dissertation, Western Michigan University, 2023.

\bibitem{Bowling1}
Bowling, A., Zonal and cozonal labelings using arbitrary abelian groups. \textit{Util. Math.}, 124:83-100, 2025.

\bibitem{BowLow}
Bowling, A. and Low, R., Cozonality in Plane Graphs of Maximum Degree Four. \textit{Discrete Math. Lett.}, 17:1-6, 2026.

\bibitem{BowXie}
Bowling, A. and Xie, W., Zonal labelings and Tait colorings from a new perspective. \textit{Aequationes Math.}, 98:1611-1625, 2024.

\bibitem{BowZhang}
Bowling, A. and Zhang, P., Absolutely and conditionally zonal graphs. \textit{Electronic J. Math.}, 4:1-11, 2022.

\bibitem{BowXieLow}
Bowling, A., Xie, W. and Low, R.M., Cozonal labelings of plane graphs. 
\textit{Bull. Inst. Combin. Appl.}, 105:54-69, 2025.

\bibitem{CharEganZhang}
Chartrand, G., Egan, C. and Zhang, P., \textit{How to Label a Graph}, Springer, New York, 2019.

\bibitem{CharEganZhangRevis}
Chartrand, G., Egan, C. and Zhang, P., Zonal graphs revisited. \textit{Bull. Inst. Combin. Appl.}, 99:133-152, 2023.

\bibitem{ChLe}
Chartrand, G. and Lesniak, L., \textit{Graphs and Digraphs}, 2nd ed., Wadsworth \& Brooks/Cole, Pacific Grove, California, 1986.

\bibitem{Dickenstein}
Dickenstein, A. and Tobis, E., Additive edge labelings. \textit{Discrete Appl. Math.}, 158:444-452, 2010.

\bibitem{DummitFoote}
Dummit, D. S. and Foote, R. M., \textit{Abstract Algebra}, 3rd ed., Wiley, Hoboken, New Jersey, 2004.

\bibitem{Grossman}
Grossman, J., Kulkarni, D. and Schochetman, I., On the minors of an incidence matrix and its Smith normal form. \textit{Linear Algebra Appl.}, 218:213-224, 1995.

\bibitem{Jeurissen}
Jeurissen, R., The incidence matrix and labelings of a graph. \textit{J. Combin. Theory Ser. B}, 30:290-301, 1981.

\bibitem{Kaibel}
Kaibel, V. and Pfetsch, M., Computing the face lattice of a polytope from its vertex-facet incidences. \textit{Comput. Geom.}, 23:281-290, 2002.

\bibitem{Yan}
Yan, W. and Li, S., On the vertex-face graphs of triangulations. \textit{Electron. J. Linear Algebra}, 36:616-628, 2020.


\end{thebibliography}
\end{document}